\documentclass[12pt]{article}
\usepackage{latexsym,amssymb,amsmath,ulem}
\usepackage{color}
\usepackage{graphicx}
\usepackage{latexsym,amssymb,amsmath,ulem}
\usepackage{amsfonts}
\usepackage{latexsym}
\usepackage{graphicx}
\usepackage{color}

\def\tto{\;{\lower 1pt \hbox{$\rightarrow$}}\kern -10pt
\hbox{\raise 2pt \hbox{$\rightarrow$}}\;}

\newtheorem{theorem}{Theorem}[section]
\newtheorem{proposition}{Proposition}[section]
\newtheorem{corollary}{Corollary}[section]
\newtheorem{lemma}{Lemma}[section]
\newtheorem{remark}{Remark}[section]
\newtheorem{example}{Example}[section]

\numberwithin{equation}{section}

\normalsize

\input{tcilatex}
\begin{document}

\title{Primal-dual optimization conditions for the robust sum of functions
with applications }
\author{N. Dinh\thanks{
International University, Vietnam National University - HCMC, Linh Trung
ward, Thu Duc district, Ho Chi Minh city, Vietnam (ndinh02@gmail.com). }, \
\ M.A. Goberna\thanks{%
Department of Mathematics, University of Alicante, Spain (mgoberna@ua.es){}}%
, \ \ M. Volle\thanks{%
Avignon University, LMA EA 2151, Avignon, France
(michel.volle@univ-avignon.fr)} }
\maketitle
\date{}

\begin{abstract}
This paper associates a dual problem to the minimization of an arbitrary
linear perturbation of the robust sum function introduced in \cite{DGV19}.
It provides an existence theorem for primal optimal solutions and, under
suitable duality assumptions, characterizations of the primal-dual optimal
set, the primal optimal set, and the dual optimal set, as well as a formula
for the subdiffential of the robust sum function. The mentioned results are
applied to get simple formulas for the robust sums of subaffine functions (a
class of functions which contains the affine ones) and to obtain conditions
guaranteeing the existence of best approximate solutions to inconsistent
convex inequality systems.
\end{abstract}

%\title{ Robust Optimization Duality }

\textperiodcentered\ Optimality conditions \textperiodcentered\ Existence of
optimal solutions \textperiodcentered\ Inconsistent convex inequality
systems \textperiodcentered\ Best approximation


\section{Introduction}

In our previous paper \cite{DGV19} we have introduced the so-called \textit{%
robust sum} $\sum\nolimits_{i\in I}^{R}{{f_{i}}}$ of an infinite family $%
\left( f_{i}\right) _{i\in I}$ of proper functions from a given locally
convex Hausdorff topological vector space $X$ to $\mathbb{R\cup }\left\{
+\infty \right\} .$ To this aim we denoted by $\mathcal{F}\left( I\right) $
the collection of all nonempty finite subsets of $I$ and defined the \textit{%
robust sum} of $\left( f_{i}\right) _{i\in I}$ as

\begin{equation*}
\sum\nolimits_{i\in I}^{R}{{f_{i}\left( x\right) :=}}\sup\limits_{J\in
\mathcal{F}\left( I\right) }\sum\nolimits_{i\in J}{{f_{i}}}\left( x\right)
,\forall x\in X.
\end{equation*}%
In order to motivate this definition, consider the finite sum $%
\sum\nolimits_{i\in J}{{f_{i}}}$ for each $J\in \mathcal{F}\left( I\right) $
and interpret $\mathcal{F}\left( I\right) $ as an \textit{uncertainty set}
for the \textit{uncertain optimization problem}%
\begin{equation*}
\mathrm{\left( P_{J}\right) }\ \ \ \ \ \ f\left( x\right) =\
\inf\limits_{x\in X}\sum\nolimits_{i\in J}{{f_{i}}}\left( x\right) .
\end{equation*}%
Then, the \textit{robust} (or \textit{pessimistic}) \textit{counterpart} of
this parametric problem is (see \cite{BEN09} and references therein) the
deterministic problem
\begin{equation}
\mathrm{\left( RP\right) }\ \ \ \ \ \ \ \inf\limits_{x\in
X}\sup\limits_{J\in \mathcal{F}\left( I\right) }\sum\nolimits_{i\in J}{{f_{i}%
}}\left( x\right) ,\text{ }  \label{RP}
\end{equation}%
whose objective function $\sum\nolimits_{i\in I}^{R}{{f_{i}}}$ cannot be
exactly computed at a given $x$ but can be approximated through the finite
sums $\sum\nolimits_{i\in J}{{f_{i}}}\left( x\right) ,$ with $J\in \mathcal{F%
}\left( I\right) .$ Observe that the above uncertain problem only makes
sense when $I$ is infinite as, otherwise, $\sum\nolimits_{i\in I}{{f_{i}}}%
\left( x\right) $ is computable at any $x\in \mathbb{R}^{n}$ and $\mathrm{%
\left( P_{I}\right) }\ $is the deterministic problem to be solved. However,
this uninteresting case allows to appreciate the pessimistic character of $%
\mathrm{\left( RP\right) }$ in comparison with $\mathrm{\left( P_{I}\right) .%
}$ Indeed, defining $I\left( x\right) :=\left\{ i\in I:f_{i}\left( x\right)
\geq 0\right\} ,$ the objective function of $\mathrm{\left( RP\right) }$
reads%
\begin{equation*}
f\left( x\right) =\left\{
\begin{array}{ll}
\max\nolimits_{i\in I}{{f_{i}}}\left( x\right) , & \text{if }I\left(
x\right) =\emptyset , \\
\sum\nolimits_{i\in I\left( x\right) }{{f_{i}}}\left( x\right) , & \text{%
else,}%
\end{array}%
\right.
\end{equation*}%
with $f$ being an upper estimate of$\ \sum\nolimits_{i\in I}{{f_{i}}}$ (the
difference $f-\sum\nolimits_{i\in I}{{f_{i}}}$ may be quite large).

It is worth observing that, in contrast with the well-known \textit{limit sum%
}
\begin{equation*}
\sum\limits_{i\in I}{{f_{i}}}\left( x\right) :=\lim\limits_{J\in \mathcal{F}%
\left( I\right) }\sum\limits_{i\in J}{{f_{i}}}\left( x\right) ,\forall x\in X
\end{equation*}%
(where $\mathcal{F}\left( I\right) $\ and $\lim $ must be interpreted as a
set directed by inclusion and the limit of the corresponding net,
respectively), the robust sum $\sum\nolimits_{i\in I}^{R}{{f_{i}}}$ is
always well-defined on $X.$

In \cite[Section 1]{DGV19} we gave two examples of optimization problems
arising in extended regression and best approximate solution to inconsistent
linear system\ which can be formulated as $\mathrm{\left( RP\right) ,}$ with
$\left( f_{i}\right) _{i\in I}$ being families of quadratic functions and
maxima of affine functions, respectively.

In this paper we assume that some element $\overline{x}^{\ast }$ of the dual
space $X^{\ast }$ of $X$ is given and introduce a dual problem for the
\textit{linearly perturbed robust sum} $\sum\nolimits_{i\in I}^{R}{{f_{i}}}%
-\left\langle \overline{x}^{\ast },\cdot \right\rangle .$ More precisely, we
are concerned with the non-emptiness and the structure of the optimal sets
of the dual pair of optimization problems%
\begin{equation*}
(\mathrm{RP}_{\overline{x}^{\ast }})\ \ \ \ \ \ \ \inf \left\{ {f}\left(
x\right) -\left\langle \overline{x}^{\ast },x\right\rangle :x\in X\right\}
\end{equation*}%
and%
\begin{equation*}
(\mathrm{RD}_{\overline{x}^{\ast }})\quad \sup \left\{ -\sum_{j\in
J}f_{j}^{\ast }(x_{j}^{\ast }):\left( J,(x_{j}^{\ast })_{j\in J}\right) \in
\mathbb{F}\left( \overline{x}^{\ast }\right) \right\} ,
\end{equation*}%
where $f:=\sum\nolimits_{i\in I}^{R}{{f_{i}}}$ represents the robust sum of
the family $\left( f_{i}\right) _{i\in I},$ the objective function $%
-\sum_{j\in J}f_{j}^{\ast }(x_{j}^{\ast })$ of $(\mathrm{RD}_{\overline{x}%
^{\ast }})$ is well defined thanks to the properness of $f_{i}$
(guaranteeing that its conjugate function $f_{i}^{\ast }$ does not take the
value $-\infty )$ for all $i\in I,$ and the feasible set of the dual
problem, $\mathbb{F}\left( \overline{x}^{\ast }\right) ,$ is defined as
\begin{equation*}
\mathbb{F}\left( \overline{x}^{\ast }\right) :=\left\{ \left( J,(x_{j}^{\ast
})_{j\in J}\right) :J\in \mathcal{F}\left( I\right) ,(x_{j}^{\ast })_{j\in
J}\in \left( X^{\ast }\right) ^{J},\sum_{j\in J}x_{j}^{\ast }=\overline{x}%
^{\ast }\right\} .
\end{equation*}

When $\overline{x}^{\ast }$ is the null functional, the pair formed by $(%
\mathrm{RP}_{\overline{x}^{\ast }})$ and $(\mathrm{RD}_{\overline{x}^{\ast
}})$ collapses to the pair of dual problems analyzed in \cite{DGV19},\ for
which we characterized weak duality, zero duality gap, and strong duality,
and their corresponding stable versions, but without paying attention to
their optimal solution sets.

Many works have been written on the numerical methods for the problem of
best least squares solutions of inconsistent finite linear inequality
systems (see, e.g., \cite{PS14} and references therein), for which the
existence of optimal solutions has been proved in three different ways in
\cite{CHP2017}. Unfortunately, as shown in \cite{GHL18}, the existence of
optimal solution for the best least squares approximation problems relies on
the finiteness of the number of constraints and the type of norm used to
measure the residual of an approximate solution. The novelties of Section 6,
in comparison with its unique antecedent \cite{GHL18}, is that, here, we
consider convex systems instead of linear ones, describe the structure of
the sets of best $\ell _{1}$ and $\ell _{\infty }$ approximate solutions
(instead of just an existence theorem for best $\ell _{\infty }$
approximation problems), and provide strong duality theorems for best $\ell
_{1}$ and $\ell _{\infty }$ approximation problems.

This paper is organized as follows. Section 2 introduces the necessary
notation and some preliminary results. Section 3 provides an existence
theorem for primal optimal solutions. Section 4 characterizes the
primal-dual optimal solutions with zero duality gap, as well as, under
suitable assumptions, primal optimal solutions, dual optimal solutions and
also provides a closed formula for the subdifferential of the robust sum
function. Section 5 provides formulas for the robust sums of subaffine
functions (concept introduced in Section 2). Finally, Section 6 provides
existence theorems for best approximate solutions to inconsistent convex
inequality systems with respect to the $\ell _{\infty }$ and the $\ell _{1}$
pseudo-norms.

\section{Preliminaries}

We first recall some standard notation regarding locally convex spaces to be
used in the sequel. We denote by $0_{X}$ and \ $0_{X}^{\ast }$ the null
vectors of $X$ and $X^{\ast },$ respectively. Given a set $A\subset X,$ we
denote by $\limfunc{co}A$, $\limfunc{cone}A,$ $\limfunc{aff}A,$ $\overline{A}%
,$ $\overline{\limfunc{co}}A,$ and $\overline{\limfunc{cone}}A$ the convex
hull of $A$, the cone generated by $A\cup \left\{ 0_{X}\right\} ,$ the
smallest linear manifold containing $A,$ the closure of $A,$ the closed
convex hull $A$, and the closed conic hull of $A,$ respectively. The same
notation is used when either $A\subset X^{\ast }$ (by default equipped
equipped with the $w^{\ast }-$topology) or $A\subset X^{\ast }\times \mathbb{%
R}$ (equipped with the product topology). We represent by $\limfunc{proj}%
\nolimits_{X^{\ast }}$ the mapping from $X^{\ast }\times \mathbb{R}$ to $%
X^{\ast }$ such that $\limfunc{proj}\nolimits_{X^{\ast }}\left( x^{\ast
},r\right) =x^{\ast }.$ When $X=\mathbb{R}^{n},$ we denote by $\limfunc{ri}A$
the relative interior of $A.$

Given $A,B\subset X,$ $A$ is said \cite{Bot10} to be \textit{closed
regarding to} $B$ if $B\cap \overline{A}=B\cap A.$ Clearly, $A$ is closed
regarding $B$ if and only if $A$ is closed regarding each subset of $B.$

We denote by $\overline{\mathbb{R}}$ the extended real line with $\pm \infty
$ and by $\overline{\mathbb{R}}^{X}$ the linear space of functions from $X$
to $\overline{\mathbb{R}}.$ Given $h\in \overline{\mathbb{R}}^{X},$ its
\textit{lower level sets} are $\left[ h\leq r\right] :=\{x\in X:h(x)\leq
r\}, $ with $r\in \mathbb{R},$\ its \textit{domain} is the set $\limfunc{dom}%
h:=\{x\in X:h(x)<+\infty \},$ its \textit{epigraph} is $\limfunc{epi}%
h:=\{\left( x,r\right) \in X\times \mathbb{R}:h(x)\leq r\},$ its \textit{%
strict epigraph} is $\limfunc{epi}\nolimits_{s}h:=\{\left( x,r\right) \in
X\times \mathbb{R}:h(x)<r\},$ and its \textit{Fenchel conjugate} the
function $h^{\ast }\in \overline{\mathbb{R}}^{X^{\ast }}$ such that
\begin{equation*}
h^{\ast }(x^{\ast }):=\sup \{\langle x^{\ast },x\rangle -h(x):x\in
X\},\forall x^{\ast }\in X^{\ast }.
\end{equation*}%
Moreover, the \textit{closed hull} of $h$ is the function $\overline{h}\in
\overline{\mathbb{R}}^{X}$ whose epigraph $\limfunc{epi}\overline{h}$ is the
closure of $\limfunc{epi}h$ in $X\times \mathbb{R}.$ The definitions are
similar if $h\in \overline{\mathbb{R}}^{X^{\ast }};$ in particular, $%
\overline{h}$ is the $w^{\ast }-$closed hull of $h.$ The \textit{%
subdifferential} of $h$ at $a\in X$ is%
\begin{equation*}
\partial h(a):=\left\{
\begin{array}{ll}
\{x^{\ast }\in X^{\ast }\,:\,h(x)\geq h(a)+\langle x^{\ast },x-a\rangle
,\forall x\in X\}, & \text{if }h(a)\in \mathbb{R}, \\
\emptyset , & \text{else.}%
\end{array}%
\right.
\end{equation*}%
The \textit{indicator function} of $A\subset X$ \ is represented by $\delta
_{A}$ (i.e. $\delta _{A}(x)=0$ if $x\in A,$ and $\delta _{A}(x)=+\infty $ if
$x\notin A$). The \textit{support function} of $A\neq \emptyset ,$ $\sigma
_{A}\left( x^{\ast }\right) :=\sup\limits_{x\in A}\langle x^{\ast },x\rangle
,$ is the conjugate of its indicator, i.e., $\sigma _{A}=$ $\delta
_{A}^{\ast }.$ The support functions are \textit{sublinear}, i.e., they are
subaditive and positively homogeneous.

We denote by $\Gamma \left( X\right) $ the cone of $\overline{\mathbb{R}}%
^{X} $ formed by the proper closed convex functions on $X.$ For instance, $%
\delta _{A}\in \Gamma \left( X\right) $ if and only if $A$ is a nonempty
closed convex set while $\sigma _{A}\in \Gamma \left( X^{\ast }\right) $ for
all nonempty $A\subset X.$ The sublinear elements of $\Gamma \left( X\right)
$ are the support functions of the nonempty $w^{\ast }-$closed convex
subsets of $X^{\ast }.$

The continuous affine functions\ on $X$ are the sums of continuous linear
functionals with constants, i.e., functions of the form $\left\langle
a^{\ast },\cdot \right\rangle +r=\sigma _{\left\{ a^{\ast }\right\} }+r,$
with $a^{\ast }\in X^{\ast }$ and $r\in \mathbb{R}.$ In the same vein, we
define the \textit{subaffine functions} on $X$ as those functions which can
be expressed as $\sigma _{A}+r,$ with $A$ being a nonempty $w^{\ast }-$%
closed convex subset of $X^{\ast }$ and $r\in \mathbb{R}.$ For instance, the
polar $A^{\circ }$ of such a set $A$ is the lower level set of some
subaffine function.\ Indeed,
\begin{equation*}
A^{\circ }:=\left\{ x\in X:\left\langle a^{\ast },x\right\rangle \leq
1,\forall a^{\ast }\in A\right\} =\left[ \sigma _{A}-1\leq 0\right] .
\end{equation*}%
Obviously, any continuous affine function is subaffine.

\begin{remark}
The above class of subaffine functions is not related with others types of
functions introduced under the same name in different settings:\newline
1. Generalized convexity (see, e.g., \cite{PV87}, \cite{MLegaz88},\cite%
{Penot98}, \cite{RD05}): a function $f\in \mathbb{R}^{X}$ is called
subaffine (or truncated affine) if it can be written as $f=\min \left\{
x^{\ast }+r,s\right\} ,$ for $x^{\ast }\in X^{\ast }$ and $r,s\in \mathbb{R}%
. $ \newline
2. Elliptic PDEs (see, e.g., \cite{HL09}, \cite{NV10}): a function $f\in
\mathbb{R}^{\mathbb{R}^{n}}$ is called subaffine if it is upper
semicontinuous and there exists a ball $B$\ such that for each affine
function $h$, $f\leq h$ on $\limfunc{bd}B$ implies that $f\leq h$ on $B$. A $%
\mathcal{C}^{2}$ function is subaffine in this sense iff its Hessian matrix
has at least one nonnegative eigenvalue at each point.
\end{remark}

We now come back to the pair of problems $(\mathrm{RP}_{\overline{x}^{\ast
}})$ and $(\mathrm{RD}_{\overline{x}^{\ast }}),$ whose optimal sets are
respectively denoted%
\begin{equation*}
\limfunc{sol}(\mathrm{RP}_{\overline{x}^{\ast }})=\left\{ x\in X:{f}\left(
x\right) -\left\langle \overline{x}^{\ast },x\right\rangle =\inf (\mathrm{RP}%
_{\overline{x}^{\ast }})\right\}
\end{equation*}%
and%
\begin{equation*}
\limfunc{sol}(\mathrm{RD}_{\overline{x}^{\ast }})=\left\{ \left(
J,(x^*_{j})_{j\in J}\right) \in \mathbb{F}\left( \overline{x}^{\ast }\right)
:-\sum_{j\in J}f_{j}^{\ast }(x_{j}^{\ast })=\sup (\mathrm{RD}_{\overline{x}%
^{\ast }})\right\} .
\end{equation*}%
When $\limfunc{sol}(\mathrm{RP}_{\overline{x}^{\ast }})\neq \emptyset $ we
write $\min (\mathrm{RP}_{\overline{x}^{\ast }})$ instead of $\inf (\mathrm{%
RP}_{\overline{x}^{\ast }}).$ Similarly, we write $\max (\mathrm{RD}_{%
\overline{x}^{\ast }})$ instead of $\sup (\mathrm{RD}_{\overline{x}^{\ast
}}) $ if $\limfunc{sol}(\mathrm{RD}_{\overline{x}^{\ast }})\neq \emptyset .$

Adopting the robust optimization approach under uncertainty (as in \cite%
{BJL13}, \cite{DGLV18}, \cite{DGLV18B}, \cite{LJL11}, etc.) we have shown in
\cite{DGV19} that $(\mathrm{RP}_{\overline{x}^{\ast }})$ may be interpreted
as the robust optimization counterpart of some uncertain optimization
problem and $(\mathrm{RD}_{\overline{x}^{\ast }})$ as its optimistic dual.
In particular, the relation%
\begin{equation}
\sup (\mathrm{RD}_{\overline{x}^{\ast }})\leq \inf (\mathrm{RP}_{\overline{x}%
^{\ast }})  \label{1.1}
\end{equation}%
always holds \cite[Proposition 3.1]{DGV19}. The characterization of the%
\textit{\ strong duality}, namely $\inf (\mathrm{RP}_{\overline{x}^{\ast
}})=\max (\mathrm{RD}_{\overline{x}^{\ast }}),$ involves the set%
\begin{equation}
\mathcal{A}:=\dbigcup\limits_{\QATOP{J\in \mathcal{F}\left( I\right) }{%
\hfill }}\sum\limits_{j\in J}\limfunc{epi}f_{j}^{\ast }.  \label{eqA}
\end{equation}%
As shown below, the set $\mathcal{A}$\ may be convex in favorable
circumstances.

\begin{lemma}
\label{Lemma1.1}Let $\left( A_{i}\right) _{i\in I}$ be a family of convex
subsets of a linear space $Z$ such that $0_{Z}\in \bigcap\limits_{i\in
I}A_{i}.$ Then $A:=\dbigcup\limits_{\QATOP{J\in \mathcal{F}\left( I\right) }{%
\hfill }}\sum\limits_{j\in J}A_{j}$ is a convex subset of $Z.$
\end{lemma}

\textit{Proof}. Notice that $\left( \sum\limits_{j\in J}A_{j}\right) _{J\in
\mathcal{F}\left( I\right) }$ is a family of convex subsets of $Z$ which is
directed with respect to the inclusion. It follows that $A$ is
convex.\noindent $\hfill \square \medskip $

\begin{example}
\label{Ex1.1}The set $\mathcal{A}=\dbigcup\limits_{\QATOP{J\in \mathcal{F}%
\left( I\right) }{\hfill }}\sum\limits_{j\in J}\limfunc{epi}f_{j}^{\ast }$
is convex if the functions $f_{j},$ $j\in J,$ are non-negative.
\end{example}

\begin{example}
\label{Ex1.2}The set $A:=\dbigcup\limits_{\QATOP{J\in \mathcal{F}\left(
I\right) }{\hfill }}\sum\limits_{j\in J}\func{dom}f_{j}^{\ast }$ is convex
if each function $f_{j},$ $j\in J,$ is bounded below.
\end{example}

We have the following characterization of strong duality under convexity.

\begin{theorem}[Strong zero duality gap under convexity]
\label{Theor1.1}\textrm{\cite[Theorem 6.1]{DGV19}} Assume the $f_{i}\in $ $%
\Gamma \left( X\right) ,$ $i\in I,$ and $\func{dom}f\neq \emptyset .$ The
next statements are equivalent:\newline
$(i)$ $\inf (\mathrm{RP}_{\overline{x}^{\ast }})=\max (\mathrm{RD}_{%
\overline{x}^{\ast }}).$ \newline
$(ii)$ $\mathcal{A}$ is $w^{\ast }-$closed convex regarding $\left\{
\overline{{x}}^{\ast }\right\} \times \mathbb{R}.$\newline
In particular, $(i)$ holds for any $\overline{{x}}^{\ast }\in X^{\ast }$ if
and only if $\mathcal{A}$ is $w^{\ast }-$closed convex.
\end{theorem}

\section{Minimizing the robust sum: existence of primal optimal solutions}

In this section we assume that $\left( f_{i}\right) _{i\in I}\subset \Gamma
\left( X\right) $ and, unless specified otherwise, that $f=\sum\nolimits_{i%
\in I}^{R}{{f_{i}}}$ is proper. We thus have $f\in \Gamma \left( X\right) .$
Additionally, we suppose that%
\begin{equation}
f\text{ is weakly inf-locally compact}  \label{2.1}
\end{equation}%
in the sense that the lower level set $\left[ f\leq r\right] $ is weakly
locally compact for each $r\in \mathbb{R}.$ Let us note that this condition
is always satisfied if $X$ is finitely dimensional. It is also satisfied if $%
\sup\nolimits_{i\in I}f_{i}$ is weakly inf-locally compact or, a fortiori,
if there exists $i\in I$ such that $f_{i}$ is weakly inf-locally compact.

By \cite[Chapter 1, Proposition 5.4]{Joly70} or by \cite[Theorem 7.7.6]%
{Lau72}, (\ref{2.1}) is equivalent to:
\begin{equation*}
f^{\ast }\text{ is quasicontinuous with respect to the Mackey topology }\tau
\left( X^{\ast },X\right) \text{ on }X^{\ast }.
\end{equation*}

Let us recall that a convex function $\xi :X^{\ast }\longrightarrow
\overline{\mathbb{R}}$ is said to be $\tau \left( X^{\ast },X\right) -$%
\textit{quasicontinuous} if the following four properties are satisfied (%
\cite{Joly70}, \cite{JL71}, \cite{Lau72}):

\begin{itemize}
\item $\limfunc{aff}\func{dom}\xi $ is $\tau \left( X^{\ast },X\right) -$%
closed (or $w^{\ast }-$closed).

\item $\limfunc{aff}\func{dom}\xi $ is of finite codimension.

\item The $\tau \left( X^{\ast },X\right) -$relative interior of $\func{dom}%
\xi ,$ say $\limfunc{ri}\func{dom}\xi ,$ is nonempty.

\item The restriction of $\xi $ to $\limfunc{aff}\func{dom}\xi $ is $\tau
\left( X^{\ast },X\right) -$continuous on $\limfunc{ri}\func{dom}\xi .$
\end{itemize}

\begin{remark}
\label{rem3.0} A convex function majorized by a $\tau (X^{\ast },X)$%
-quasicontinuous one is\newline
$\tau (X^{\ast },X)$-quasicontinuous, too (see \cite[Theorem 2.4]{MV97},
\cite[Proposition 2.2.15]{Za02}). If $X=X^{\ast }=\mathbb{R}^{n}$, any
extended real-valued convex function with nonempty domain is quasicontinuous.
\end{remark}

Let us consider the subdifferential of $f^{\ast }$ at $\overline{x}^{\ast
}\in X^{\ast },$ namely,%
\begin{equation*}
\partial f^{\ast }\left( \overline{x}^{\ast }\right) =\left\{
\begin{array}{ll}
\left\{ x\in X:f^{\ast }\left( x^{\ast }\right) \geq f^{\ast }\left(
\overline{x}^{\ast }\right) +\left\langle x^{\ast }-\overline{x}^{\ast
},x\right\rangle ,\forall x^{\ast }\in X^{\ast }\right\} , & \text{if }%
f^{\ast }\left( \overline{x}^{\ast }\right) \in \mathbb{R}, \\
\emptyset , & \text{else.}%
\end{array}%
\right.
\end{equation*}

For $\overline{x}^{\ast }\in \func{dom}f^{\ast },$ since $f\in \Gamma \left(
X\right) $ entails $f^{\ast \ast }=f,$ one has%
\begin{equation}
\partial f^{\ast }\left( \overline{x}^{\ast }\right) =\func{argmin}\left(
f-\left\langle \overline{x}^{\ast },\cdot \right\rangle \right) =\limfunc{sol%
}(\mathrm{RP}_{\bar{x}^{\ast }}).  \label{solRPx*}
\end{equation}

We are faced with the subdifferentiability of $f^{\ast }$ at $\overline{x}%
^{\ast },$ for which the dual version \cite[Theorem III.3]{MV97} gives a
very useful criterion:

\begin{lemma}
\label{Lemma2.1}Assume that $g\in \Gamma \left( X\right) $ is weakly
inf-locally compact and%
\begin{equation}
\overline{\limfunc{cone}}\left( \func{dom}g^{\ast }-\overline{x}^{\ast
}\right) \text{ is a linear subspace of }X^{\ast }.  \label{2.2}
\end{equation}%
Then $\partial g^{\ast }\left( \overline{x}^{\ast }\right) $ is the sum of a
nonempty weakly compact convex set and a finitely dimensional linear
subspace of $X.$
\end{lemma}

\begin{remark}
\label{rem3.2b} Condition (\ref{2.2}) means that the sets $\limfunc{dom}%
g^{\ast }$ and $\{x^{\ast }\}$ are united in the sense that they cannot be
properly separated (all weak$^{\ast }$-closed hyperplanes which separate
them contain both of them). A sufficient (in general not necessary)
condition for this is that $x^{\ast }$ belongs to the relative algebraic
interior of $\limfunc{dom}g^{\ast }$ (see \cite[Proposition 1.2.8]{Za02} for
more details).
\end{remark}

To exploit Lemma \ref{Lemma2.1} in the case that $g=f=\sum\nolimits_{i\in
I}^{R}{{f_{i},}}$ we need an explicit formulation of the criterion (\ref{2.2}%
) in terms of the functions $f_{i}^{\ast }.$ To this end, let us consider
the function $\varphi $ defined on $X^{\ast }$ by%
\begin{equation}
\varphi \left( x^{\ast }\right) :=\inf \left\{ \sum_{j\in J}f_{j}^{\ast
}(x_{j}^{\ast }):\left( J,(x_{j}^{\ast })_{j\in J}\right) \in \mathbb{F}%
\left( x^{\ast }+\overline{x}^{\ast }\right) \right\} ,\forall x^{\ast }\in
X^{\ast }.  \label{2.0}
\end{equation}%
One has straightfordwardly
\begin{equation*}
\varphi ^{\ast }\left( x\right) =f\left( x\right) -\left\langle \overline{x}%
^{\ast },x\right\rangle ,\forall x\in X,
\end{equation*}%
\begin{equation*}
\varphi ^{\ast \ast }\left( x^{\ast }\right) =f^{\ast }\left( x^{\ast }+%
\overline{x}^{\ast }\right) ,\forall x^{\ast }\in X^{\ast },
\end{equation*}%
and
\begin{equation}
\func{dom}f^{\ast }-\overline{x}^{\ast }=\func{dom}\varphi ^{\ast \ast }.
\label{domf*}
\end{equation}

Since $\func{dom}\varphi ^{\ast }=\func{dom}f\neq \emptyset ,$ the
biconjugate function $\varphi ^{\ast \ast }$ coincides with the $w^{\ast }-$%
closed convex hull $\overline{\limfunc{co}}\varphi $ of $\varphi ,$ which
satisfies%
\begin{equation}
\limfunc{epi}\overline{\limfunc{co}}\varphi =\overline{\limfunc{co}}\limfunc{%
epi}\varphi .  \label{2.3}
\end{equation}

Let us observe that
\begin{equation}
\limfunc{proj}\nolimits_{X^{\ast }}\left( \limfunc{co}\limfunc{epi}\varphi
\right) =\limfunc{co}\func{dom}\varphi .  \label{2.4}
\end{equation}

Now, by (\ref{2.3}) and (\ref{2.4}), one has%
\begin{equation*}
\func{dom}\overline{\limfunc{co}}\varphi =\limfunc{proj}\nolimits_{X^{\ast
}}\left( \overline{\limfunc{co}}\limfunc{epi}\varphi \right) \subset
\overline{\limfunc{proj}\nolimits_{X^{\ast }}\left( \limfunc{co}\limfunc{epi}%
\varphi \right) }=\overline{\limfunc{co}}\func{dom}\varphi ,
\end{equation*}%
and, since $\overline{\limfunc{co}}\func{dom}\varphi $ is $w^{\ast }-$closed,%
\begin{equation*}
\overline{\func{dom}\overline{\limfunc{co}}\varphi }\subset \overline{%
\limfunc{co}}\func{dom}\varphi .
\end{equation*}%
Conversely, since $\overline{\limfunc{co}}\varphi \leq \varphi ,$ we have $%
\func{dom}\varphi \subset \func{dom}\overline{\limfunc{co}}\varphi $ and,
since $\func{dom}\overline{\limfunc{co}}\varphi $ is convex, $\limfunc{co}%
\func{dom}\varphi \subset \func{dom}\overline{\limfunc{co}}\varphi .$ So, $%
\overline{\limfunc{co}}\func{dom}\varphi =\overline{\limfunc{co}\func{dom}%
\varphi }\subset \overline{\func{dom}\overline{\limfunc{co}}\varphi }.$
Consequently,
\begin{equation}
\overline{\limfunc{co}}\func{dom}\varphi =\overline{\func{dom}\overline{%
\limfunc{co}}\varphi },  \label{2.5}
\end{equation}%
and hence, it follows from \eqref{domf*} that
\begin{equation*}
\begin{array}{ll}
\overline{\limfunc{cone}}\left( \func{dom}f^{\ast }-\overline{x}^{\ast
}\right) & =\overline{\limfunc{cone}}\func{dom}\varphi ^{\ast \ast }=%
\overline{\limfunc{cone}}\func{dom}\overline{\limfunc{co}}\varphi \\
& =\overline{\limfunc{cone}}\left( \overline{\func{dom}\overline{\limfunc{co}%
}\varphi }\right) =\overline{\limfunc{cone}}\left( \overline{\limfunc{co}}%
\func{dom}\varphi \right) \\
& =\overline{\limfunc{cone}}\left( \limfunc{co}\func{dom}\varphi \right) .%
\end{array}%
\end{equation*}%
Now, from the very definition of $\varphi ,$ one has
\begin{equation*}
\func{dom}\varphi =\left( \dbigcup\limits_{\QATOP{J\in \mathcal{F}\left(
I\right) }{\hfill }}\sum\limits_{j\in J}\func{dom}f_{j}^{\ast }\right) -%
\overline{x}^{\ast },
\end{equation*}%
and the criterion (\ref{2.2}) writes, for $g=f,$%
\begin{equation}
\overline{\limfunc{cone}}\limfunc{co}\left\{ \left( \dbigcup\limits_{\QATOP{%
J\in \mathcal{F}\left( I\right) }{\hfill }}\sum\limits_{j\in J}\func{dom}%
f_{j}^{\ast }\right) -\overline{x}^{\ast }\right\} \text{ is a linear
subspace of }X^{\ast }.  \label{2.6}
\end{equation}%
Together with \eqref{solRPx*} and Lemma \ref{Lemma2.1}, we have thus proved
the following result:

\begin{theorem}[Existence of optimal solution]
\label{Theor2.1} Assume that $\left( f_{i}\right) _{i\in I}\subset \Gamma
\left( X\right) ,$ $f=\sum\nolimits_{i\in I}^{R}{{f_{i}}}$ is proper weakly
inf-locally compact and (\ref{2.6}) holds. Then $(\mathrm{RP}_{\overline{x}%
^{\ast }})$ admits {\ at least an} optimal solution. More precisely, $%
\limfunc{sol}(\mathrm{RP}_{\overline{x}^{\ast }})$ is the sum of a nonempty
convex weakly compact set and a finitely dimensional linear subspace of $X.$
\end{theorem}

For nonnegative functions we obtain:

\begin{corollary}
\label{Corol3.1b} Let $\left( f_{i}\right) _{i\in I}$ be a family of
nonnegative $\Gamma (X)$-functions such that the infinite sum $%
\sum\nolimits_{i\in I}{{f_{i}}}$ is proper weakly inf-locally compact.
Assume that
\begin{equation}
\overline{\mathrm{cone}}\dbigcup\limits_{\QATOP{J\in \mathcal{F}\left(
I\right) }{\hfill }}\sum\limits_{j\in J}\func{dom}f_{j}^{\ast }\ \ \text{ is
a linear subspace of }X^{\ast }.  \label{3.10b}
\end{equation}%
Then the optimal solution set of the problem
\begin{equation*}
\inf_{x\in X}\sum\nolimits_{i\in I}f_{i}(x)
\end{equation*}%
is the sum of a nonempty convex weakly compact set and a finitely
dimensional linear subspace of $X$.
\end{corollary}

\textit{Proof}. Since the functions $f_{i}$, $i\in I$ are nonnegative, their
robust sum coincides with the infinite sum $\sum\nolimits_{i\in I}{{f_{i}}}$%
. Moreover, one has $0_{X^{\ast }}\in \limfunc{dom}f_{i}^{\ast }$ for each $%
i\in I$, and the set $\dbigcup\limits_{\QATOP{J\in \mathcal{F}\left(
I\right) }{\hfill }}\sum\limits_{j\in J}\func{dom}f_{j}^{\ast }$ is convex
(see Example \ref{Ex1.2}). We conclude the proof with Theorem \ref{Theor2.1}%
. $\hfill \square \medskip $

\begin{remark}
\label{rem33b} If $I$ is finite and all functions $f_{i},$ $i\in I,$ are
nonnegative, then
\begin{equation*}
\dbigcup\limits_{\QATOP{J\in \mathcal{F}\left( I\right) }{\hfill }%
}\sum\limits_{j\in J}\func{dom}f_{j}^{\ast }=\sum\limits_{i\in I}\limfunc{dom%
}f_{i}^{\ast },
\end{equation*}%
and condition \eqref{3.10b} becomes
\begin{equation*}
\overline{\limfunc{cone}}\sum\limits_{i\in I}\limfunc{dom}f_{i}^{\ast }\
\text{ is a linear subspace of }X^{\ast }.
\end{equation*}
\end{remark}

Observe that, under the assumptions of Theorem \ref{Theor2.1}, one has in
particular $\inf (\mathrm{RP}_{\overline{x}^{\ast }})\in \mathbb{R}.$
Observe also that when $X=X^{\ast }=\mathbb{R}^{n},$ (\ref{2.2}) writes $%
\overline{x}^{\ast }\in \limfunc{ri}\left( \func{dom}g^{\ast }\right) ,$ and
in such a case, one has the next corollary.

\begin{corollary}
\label{Corol2.1} Assume that $\left( f_{i}\right) _{i\in I}\subset \Gamma
\left( \mathbb{R}^{n}\right) ,$ $\func{dom}f\neq \emptyset ,$ and
\begin{equation}
\overline{x}^{\ast }\in \limfunc{ri}\limfunc{co}\left( \dbigcup\limits_{%
\QATOP{J\in \mathcal{F}\left( I\right) }{\hfill }}\sum\limits_{j\in J}\func{%
dom}f_{j}^{\ast }\right) .  \label{2.7}
\end{equation}%
Then, $\limfunc{sol}(\mathrm{RP}_{\overline{x}^{\ast }})$ is the sum of a
nonempty convex compact set and a linear subspace of $\mathbb{R}^{n}.$
\end{corollary}

\begin{remark}
If each function $f_{i},$ $i\in I,$ is bounded below, then (see Example 2.2)
the criteria (\ref{2.6}) and (\ref{2.7}) collapse respectively to%
\begin{equation*}
\overline{\limfunc{cone}}\left\{ \left( \dbigcup\limits_{\QATOP{J\in
\mathcal{F}\left( I\right) }{\hfill }}\sum\limits_{j\in J}\func{dom}%
f_{j}^{\ast }\right) -\overline{x}^{\ast }\right\} \text{ is a linear
subspace of }X^{\ast }
\end{equation*}%
and
\begin{equation*}
\overline{x}^{\ast }\in \limfunc{ri}\left( \dbigcup\limits_{\QATOP{J\in
\mathcal{F}\left( I\right) }{\hfill }}\sum\limits_{j\in J}\func{dom}%
f_{j}^{\ast }\right) .
\end{equation*}
\end{remark}

Note that the conclusion of Theorem \ref{Theor2.1} does not entail that%
\begin{equation}
\min (\mathrm{RP}_{\overline{x}^{\ast }})=\sup (\mathrm{RD}_{\overline{x}%
^{\ast }}).  \label{2.8}
\end{equation}%
One has in fact, with $\varphi $ defined as in (\ref{2.0}), the following
lemma.

\begin{lemma}
\label{Lemma2.2} Assume that either $\sup (\mathrm{RD}_{\overline{x}^{\ast
}})=+\infty $ or $\varphi $ is subdifferentiable at $0_{X^{\ast }}.$ Then (%
\ref{2.8}) holds.
\end{lemma}

\textit{Proof}. Since $\inf (\mathrm{RP}_{\overline{x}^{\ast }})\geq \sup (%
\mathrm{RD}_{\overline{x}^{\ast }}),$ (\ref{2.8})\ is obvious if $\sup (%
\mathrm{RD}_{\overline{x}^{\ast }})=+\infty .$ Assume now that $\overline{x}%
\in \partial \varphi \left( 0_{X^{\ast }}\right) .$ Then $\varphi \left(
0_{X^{\ast }}\right)+ \varphi ^{\ast }\left( \overline{x}\right) = \langle
0_{X^*}, \overline{x}\rangle = 0$ and we thus have%
\begin{equation*}
\inf (\mathrm{RP}_{\overline{x}^{\ast }})\leq f\left( \overline{x}\right)
-\left\langle \overline{x}^{\ast },\overline{x}\right\rangle =\varphi ^{\ast
}\left( \overline{x}\right) =-\varphi \left( 0_{X^{\ast }}\right) =\sup (%
\mathrm{RD}_{\overline{x}^{\ast }})\leq \inf (\mathrm{RP}_{\overline{x}%
^{\ast }}),
\end{equation*}%
and (\ref{2.8}) follows.$\hfill \square \medskip $

\begin{remark}
\label{rem32a} Recall that $\mathcal{A}=\dbigcup\limits_{\QATOP{J\in
\mathcal{F}\left( I\right) }{\hfill }}\sum\limits_{j\in J}\limfunc{epi}%
f_{j}^{\ast }$ and $\func{dom}\varphi =\left( \dbigcup\limits_{\QATOP{J\in
\mathcal{F}\left( I\right) }{\hfill }}\sum\limits_{j\in J}\func{dom}%
f_{j}^{\ast }\right) -\overline{x}^{\ast }.$ From (\ref{2.0}) one has%
\begin{equation*}
\limfunc{epi}\nolimits_{s}\varphi \subset \mathcal{A}-\left( \overline{x}%
^{\ast },0\right) \subset \limfunc{epi}\varphi
\end{equation*}%
and, consequently,
\begin{equation*}
\varphi \left( \overline{x}^*\right) =\inf \left\{ t\in \mathbb{R}:\left(
x^{\ast },t\right) \in \mathcal{A}-\left( \overline{x}^{\ast },0\right)
\right\} .
\end{equation*}%
It follows that, if $\mathcal{A}$ is convex, then $\varphi $ is convex too.
\end{remark}

\begin{theorem}[Primal attainment]
\label{Theor2.2} Assume that $\left( f_{i}\right) _{i\in I}\subset \Gamma
\left( X\right) $, $\varphi $ defined by \eqref{2.0} is convex and
Mackey-quasicontinuous, and that
\begin{equation}  \label{3.12nw}
\overline{\limfunc{cone}}\left\{ \left( \dbigcup\limits_{\QATOP{J\in
\mathcal{F}\left( I\right) }{\hfill }}\sum\nolimits_{j\in J}\func{dom}%
f_{j}^{\ast }\right) -\overline{x}^{\ast }\right\} \ \text{ is\ a \ linear\
subspace\ of } \ X^{\ast }.
\end{equation}
Then,
\begin{equation*}
\min (\mathrm{RP}_{\overline{x}^{\ast }})=\sup (\mathrm{RD}_{\overline{x}%
^{\ast }}).
\end{equation*}
\end{theorem}

\textit{Proof}. By Lemma \ref{Lemma2.2} one may assume that $\varphi \left(
0_{X^{\ast }}\right) \neq -\infty .$ By \cite[Theorem 3.3]{MV97} we have $%
\partial \varphi \left( 0_{X^{\ast }}\right) \neq \emptyset $ and by Lemma %
\ref{Lemma2.2} again we are done.$\hfill \square \medskip $

\begin{remark}
\label{rem32} Since for each $\left( i,x^{\ast }\right) \in I\times X^{\ast
} $ one has $\varphi \left( x^{\ast }\right) \leq f_{i}^{\ast }\left(
x^{\ast } +\overline{x}^{\ast } \right) ,$ the function $\varphi $ (assumed
to be convex) is Mackey-quasicontinuous whenever there exists $i_{0}\in I$
such that $f_{i_{0}}$ is weakly inf-locally compact (see Remark \ref{rem3.0}%
).
\end{remark}

\begin{corollary}
\label{Corol2.2}Let $\left( f_{i}\right) _{i\in I}\subset \Gamma \left(
\mathbb{R}^{n}\right) $ be such that $\dbigcup\limits_{\QATOP{J\in \mathcal{F%
}\left( I\right) }{\hfill }}\sum\limits_{j\in J}\limfunc{epi}f_{j}^{\ast }$
is convex and
\begin{equation}  \label{3.13nw}
\overline{x}^{\ast }\in \limfunc{ri}\left( \dbigcup\limits_{\QATOP{J\in
\mathcal{F}\left( I\right) }{\hfill }}\sum\limits_{j\in J}\func{dom}%
f_{j}^{\ast }\right) .
\end{equation}%
Then $\min (\mathrm{RP}_{\overline{x}^{\ast }})=\sup (\mathrm{RD}_{\overline{%
x}^{\ast }}).$
\end{corollary}

\textit{Proof}. As $\mathcal{A}=\dbigcup\limits_{\QATOP{J\in \mathcal{F}%
\left( I\right) }{\hfill }}\sum\limits_{j\in J}\limfunc{epi}f_{j}^{\ast }$
is convex, $\varphi $ is convex, too (Remark \ref{rem32a}). Moreover, as $X=%
\mathbb{R}^{n}$ and $\limfunc{dom}\varphi \neq \emptyset $, $\varphi $ is
Mackey-quasicontinuous. Now, again, as $X=\mathbb{R}^{n}$, \eqref{3.13nw} $%
\Leftrightarrow $ \eqref{3.12nw}, and the conclusion follows from Theorem %
\ref{Theor2.2}. $\hfill \square \medskip $

\section{Primal-dual optimality relations}

We need to introduce some additional notations. Given $g:X\longrightarrow
\overline{\mathbb{R}},$ we denote by $M_{g}:X^{\ast }\rightrightarrows X$
the set-valued mapping defined, for each $x^{\ast }\in X^{\ast }$, as%
\begin{equation*}
\left( M_{g}\right) \left( x^{\ast }\right) =\left\{
\begin{array}{ll}
\func{argmin}\left( g-\left\langle x^{\ast },\cdot \right\rangle \right) , &
\text{if }g^{\ast }\left( x^{\ast }\right) \in \mathbb{R}, \\
\emptyset , & \text{else.}%
\end{array}%
\right.
\end{equation*}%
In fact, $M_{g}$ is nothing else than the inverse of the subdifferential
mapping $\partial g:X\rightrightarrows X^{\ast },$ i.e.,%
\begin{equation*}
x\in \left( M_{g}\right) \left( x^{\ast }\right) \Longleftrightarrow x^{\ast
}\in \partial g\left( x\right) .
\end{equation*}%
One has $\left( M_{g}\right) \left( x^{\ast }\right) \subset \partial
g^{\ast }\left( x^{\ast }\right) $ and equality holds whenever $g=g^{\ast
\ast }$ (e.g., when $g\in \Gamma \left( X\right) $).

Given $x\in X,$ we denote by $S_{f}\left( x\right) $ the (possibly empty)
set of those $J\in \mathcal{F}\left( I\right) $ that realize the supremum in
the definition of the robust sum when $f\left( x\right) $ is finite:%
\begin{equation*}
S_{f}\left( x\right) =\left\{
\begin{array}{ll}
\left\{ J\in \mathcal{F}\left( I\right) :\sum\nolimits_{j\in J}f_{j}\left(
x\right) =f\left( x\right) \right\} , & \text{if }x\in \func{dom}f, \\
\emptyset , & \text{else.}%
\end{array}%
\right.
\end{equation*}

The inverse of the set-valued mapping $S_{f}:X\rightrightarrows \mathcal{F}%
\left( I\right) $ is denoted by $T_{f}.$ One has $T_{f}:\mathcal{F}\left(
I\right) \rightrightarrows X$ and
\begin{equation*}
x\in T_{f}\left( J\right) \Longleftrightarrow J\in S_{f}\left( x\right) .
\end{equation*}%
If $I$ is finite one has of course $S_{f}\left( x\right) \neq \emptyset $
for each $x\in \func{dom}f.$ We now make explicit $S_{f}\left( x\right) $ in
different situations. To this aim, we introduce the supremum function $%
f_{0}:=\sup\nolimits_{i\in I}f_{i}.$

\begin{itemize}
\item If $f_{0}\left( x\right) \leq 0$ we have $f\left( x\right)
=f_{0}\left( x\right) $ \cite[Lemma 2.5]{DGV19}. Then%
\begin{equation*}
S_{f}\left( x\right) =\left\{
\begin{array}{ll}
\left\{ \left\{ j\right\} :j\in I,f_{j}\left( x\right) =f_{0}\left( x\right)
\right\} , & \text{ if }f_{0}\left( x\right) <0, \\
\left\{ J\in \mathcal{F}\left( I\right) :f_{j}\left( x\right) =0,\forall
j\in J\right\} , & \text{ if }f_{0}\left( x\right) =0.%
\end{array}%
\right.
\end{equation*}

\item If $f_{0}\left( x\right) \in \left] 0,+\infty \right[ $ we have $%
f\left( x\right) =\sum\limits_{i\in I}f_{i}^{+}\left( x\right)
:=\sum\limits_{i\in I}\max \left\{ f_{i}\left( x\right) ,0\right\} $ \cite[%
Lemma 2.5]{DGV19} and%
\begin{equation*}
S_{f}\left( x\right) =\left\{
\begin{array}{ll}
\left\{ i\in I:f_{i}\left( x\right) >0\right\} , & \text{ if this set is
finite,} \\
\emptyset , & \text{ else.}%
\end{array}%
\right.
\end{equation*}
\end{itemize}

\begin{theorem}[Primal-dual optimality with zero duality gap]
\label{Theor3.1}Assume that all functions $f_{i}$ are proper and let $x\in
\func{dom}f$ and $\left( J,(x_{j}^{\ast })_{j\in J}\right) \in \mathbb{F}%
\left( \overline{x}^{\ast }\right) .$ Next statements are equivalent:\newline
$(i)$ $x\in \limfunc{sol}(\mathrm{RP}_{\overline{x}^{\ast }}),$ $\left(
J,(x_{j}^{\ast })_{j\in J}\right) \in \limfunc{sol}(\mathrm{RD}_{\overline{x}%
^{\ast }}),$ and $\inf (\mathrm{RP}_{\overline{x}^{\ast }})=\sup (\mathrm{RD}%
_{\overline{x}^{\ast }}).$\newline
$(ii)$ $J\in S_{f}\left( x\right) $ and $x_{j}^{\ast }\in \partial
f_{j}\left( x\right) $ for all $j\in J.$\newline
$(iii)$ $x\in T_{f}\left( J\right) \cap \left( \bigcap\limits_{j\in
J}M_{f_{j}}(x_{j}^{\ast })\right) .$\newline
If $\left( f_{i}\right) _{i\in I}\subset \Gamma \left( X\right) $ we can add
\newline
$(iv)$ $x\in T_{f}\left( J\right) \cap \left( \bigcap\limits_{j\in
J}\partial f_{j}^{\ast }(x_{j}^{\ast })\right) .$\newline
\end{theorem}

\textit{Proof. \ } From the definitions of the set-valued mappings $S_{f},$ $%
T_{f},$ and $M_{f_{j}}$ it is clear that $(ii)\Longleftrightarrow (iii)$ and
$(iii)\Longleftrightarrow (iv)$ under the assumption that $\left(
f_{i}\right) _{i\in I}\subset \Gamma \left( X\right) .$

$\left[ (i)\Longrightarrow (ii)\right] $ Since $x\in \func{dom}f$ we have $%
\sum\limits_{j\in J}{{f_{j}}}\left( x\right) \in \mathbb{R}$ and
\begin{equation}
\sum\limits_{j\in J}{{f_{j}}}\left( x\right) -\left\langle \overline{x}%
^{\ast },x\right\rangle \leq f\left( x\right) -\left\langle \overline{x}%
^{\ast },x\right\rangle =\inf (\mathrm{RP}_{\overline{x}^{\ast }})=\sup (%
\mathrm{RD}_{\overline{x}^{\ast }})=-\sum_{j\in J}f_{j}^{\ast }(x_{j}^{\ast
}).  \label{3.1}
\end{equation}

By Fenchel and Young inequality we have%
\begin{equation}
-\sum_{j\in J}f_{j}^{\ast }(x_{j}^{\ast })\leq \sum\limits_{j\in J}{{f_{j}}}%
\left( x\right) -\left\langle \overline{x}^{\ast },x\right\rangle .
\label{3.2}
\end{equation}

Since $\left( J,(x_{j}^{\ast })_{j\in J}\right) \in \mathbb{F}\left(
\overline{x}^{\ast }\right) $ we have%
\begin{equation}
\sum\limits_{j\in J}\left( {{f_{j}}}\left( x\right) -\left\langle
x_{j}^{\ast },x\right\rangle \right) =\sum\limits_{j\in J}{{f_{j}}}\left(
x\right) -\left\langle \overline{x}^{\ast },x\right\rangle .  \label{3.3}
\end{equation}

Combining (\ref{3.1}), (\ref{3.2}), and (\ref{3.3}), we obtain $%
\sum\limits_{j\in J}{{f_{j}}}\left( x\right) =f\left( x\right) ,$ that means
$J\in S_{f}\left( x\right) $ and
\begin{equation*}
\sum\limits_{j\in J}\left( {{f_{j}}}\left( x\right) +f_{j}^{\ast
}(x_{j}^{\ast })-\left\langle x_{j}^{\ast },x\right\rangle \right) =0.
\end{equation*}%
By Fenchel and Young inequality all terms of the above sum are nonnegative,
hence equal to zero, that means $x_{j}^{\ast }\in \partial f_{j}\left(
x\right) $ for all $j\in J.$\newline

$\left[ (ii)\Longrightarrow (i)\right] $ Since $J\in S_{f}\left( x\right) ,$
$\sum_{j\in J}x_{j}^{\ast }=\overline{x}^{\ast },$ $x_{j}^{\ast }\in
\partial f_{j}\left( x\right) $ for all $j\in J,$ and $\left( J,(x_{j}^{\ast
})_{j\in J}\right) \in \mathbb{F}\left( \overline{x}^{\ast }\right) ,$ we
have%
\begin{equation*}
\begin{array}{ll}
{f}\left( x\right) -\left\langle \overline{x}^{\ast },x\right\rangle &
=\sum\limits_{j\in J}{{f_{j}}}\left( x\right) -\left\langle \overline{x}%
^{\ast },x\right\rangle \\
& =\sum\limits_{j\in J}\left( {{f_{j}}}\left( x\right) -\left\langle
x_{j}^{\ast },x\right\rangle \right) \\
& =-\sum\limits_{j\in J}f_{j}^{\ast }(x_{j}^{\ast }) \\
& \leq \sup (\mathrm{RD}_{\overline{x}^{\ast }}) \\
& \leq \inf (\mathrm{RP}_{\overline{x}^{\ast }}) \\
& \leq {f}\left( x\right) -\left\langle \overline{x}^{\ast },x\right\rangle .%
\end{array}%
\end{equation*}%
All terms of the above chain of inequalities are thus equal and this proves
that $(i)$ holds.\noindent $\hfill \square \medskip $

Next corollary assumes that $\inf (\mathrm{RP}_{\overline{x}^{\ast }})=\max (%
\mathrm{RD}_{\overline{x}^{\ast }})$ (i.e., strong duality), which is
characterized (in the convex case) in Theorem \ref{Theor1.1}.

\begin{corollary}
\label{Corol3.1}Assume that all functions $f_{i}$ are proper and let $x\in
\func{dom}f$ and $\inf (\mathrm{RP}_{\overline{x}^{\ast }})=\max (\mathrm{RD}%
_{\overline{x}^{\ast }}).$ Next statements are equivalent:\newline
$(i)$ $x\in \limfunc{sol}(\mathrm{RP}_{\overline{x}^{\ast }}).$\newline
$(ii)$ For all $\left( J,(x_{j}^{\ast })_{j\in J}\right) \in \limfunc{sol}(%
\mathrm{RD}_{\overline{x}^{\ast }})$ one has $J\in S_{f}\left( x\right) $
and $x_{j}^{\ast }\in \partial f_{j}\left( x\right) $ for all $j\in J.$%
\newline
$(iii)$ There exists $\left( J,(x_{j}^{\ast })_{j\in J}\right) \in \limfunc{%
sol}(\mathrm{RD}_{\overline{x}^{\ast }})$ such that $J\in S_{f}\left(
x\right) $ and $x_{j}^{\ast }\in \partial f_{j}\left( x\right) $ for all $%
j\in J.$\newline
$(iv)$ There exists $\left( J,(x_{j}^{\ast })_{j\in J}\right) \in \mathbb{F}%
\left( \overline{x}^{\ast }\right) $ such that $J\in S_{f}\left( x\right) $
and $x_{j}^{\ast }\in \partial f_{j}\left( x\right) $ for all $j\in J.$%
\newline
Moreover, for any $\left( J,(x_{j}^{\ast })_{j\in J}\right) \in \limfunc{sol}%
(\mathrm{RD}_{\overline{x}^{\ast }})$ one has%
\begin{equation}
\limfunc{sol}(\mathrm{RP}_{\overline{x}^{\ast }}) =T_{f}\left( J\right) \cap
\left( \bigcap\limits_{j\in J}M_{f_{j}}(x_{j}^{\ast })\right) .  \label{3.4}
\end{equation}
\end{corollary}

\textit{Proof. \ } $\left[ (i)\Longrightarrow (ii)\right] $ It follows from
the statement $\left[ (i)\Longrightarrow (ii)\right] $ in Theorem \ref%
{Theor3.1}.

$\left[ (ii)\Longrightarrow (iii)\right] $ It is obvious as $\limfunc{sol}(%
\mathrm{RD}_{\overline{x}^{\ast }})\neq \emptyset .$

$\left[ (iii)\Longrightarrow (iv)\right] $ It is obvious.

$\left[ (iv)\Longrightarrow (i)\right] $ Since $J\in S_{f}\left( x\right) ,$
$\sum_{j\in J}x_{j}^{\ast }=\overline{x}^{\ast },$ and $x\in
M_{f_{j}}(x_{j}^{\ast })$ for each $j\in J,$%
\begin{equation*}
\begin{array}{ll}
\inf (\mathrm{RP}_{\overline{x}^{\ast }}) & \leq f\left( x\right)
-\left\langle \overline{x}^{\ast },x\right\rangle =\sum\limits_{j\in J}{{%
f_{j}}}\left( x\right) -\left\langle \overline{x}^{\ast },x\right\rangle \\
& =\sum\limits_{j\in J}\left( {{f_{j}}}\left( x\right) -\left\langle
x_{j}^{\ast },x\right\rangle \right) \\
& =-\sum\limits_{j\in J}f_{j}^{\ast }(x_{j}^{\ast }) \\
& \leq \sup (\mathrm{RD}_{\overline{x}^{\ast }})\leq \inf (\mathrm{RP}_{%
\overline{x}^{\ast }}).%
\end{array}%
\end{equation*}%
This ensures that $f\left( x\right) -\left\langle \overline{x}^{\ast
},x\right\rangle =$ $\inf (\mathrm{RP}_{\overline{x}^{\ast }})$ and $(i)$
holds.

Let us prove the last assertion of Corollary \ref{Corol3.1}. Let $\left(
J,(x_{j}^{\ast })_{j\in J}\right) \in \limfunc{sol}(\mathrm{RD}_{\overline{x}%
^{\ast }}).$ From $\left[ (i)\Longleftrightarrow (ii)\right] $ one has $x\in
\limfunc{sol}(\mathrm{RP}_{\overline{x}^{\ast }})$ if and only if $J\in
S_{f}\left( x\right) $ and $x_{j}^{\ast }\in \partial f_{j}\left( x\right) $
for all $j\in J$ or, equivalently,
\begin{equation*}
\ \ \ \ \ \ \ \ \ \ \ \ \ \ \ \ \ \ \ \ \ \ x\in T_{f}\left( J\right) \cap
\left( \bigcap\limits_{j\in J}M_{f_{j}}(x_{j}^{\ast })\right). \hfill \ \ \
\ \ \ \ \ \ \ \ \ \ \ \ \ \ \ \ \ \ \ \ \ \ \ \ \ \ \ \ \ \ \ \ \ \ \ \ \ \
\ \ \ \ \ \square
\end{equation*}

Notice that, if $\left( f_{i}\right) _{i\in I}\subset \Gamma \left(
X\right), $ then $M_{f_{j}}(x_{j}^{\ast }) = \partial f^\ast_{j}(x_{j}^{\ast
})$ for each $j \in J$ and the equation (\ref{3.4}) writes%
\begin{equation*}
\limfunc{sol}(\mathrm{RP}_{\overline{x}^{\ast }}) = T_{f}\left( J\right)
\cap \left( \bigcap\limits_{j\in J} \partial f^\ast_{j}(x_{j}^{\ast
})\right).
\end{equation*}

\begin{corollary}
\label{Corol3.2}Assume that all functions $f_{i}$ and $f$\ are proper and
let $\left( J,(x_{j}^{\ast })_{j\in J}\right) \in \mathbb{F}\left( \overline{%
x}^{\ast }\right) $ and $\min (\mathrm{RP}_{\overline{x}^{\ast }})=\sup (%
\mathrm{RD}_{\overline{x}^{\ast }}).$ Next statements are equivalent:\newline
$(i)$ $\left( J,(x_{j}^{\ast })_{j\in J}\right) \in \limfunc{sol}(\mathrm{RD}%
_{\overline{x}^{\ast }}), $\newline
$(ii)$ For all $x\in \limfunc{sol}(\mathrm{RP}_{\overline{x}^{\ast }})$ one
has $J\in S_{f}\left( x\right) $ and $x_{j}^{\ast }\in \partial f_{j}\left(
x\right) $ for all $j\in J, $\newline
$(iii)$ There exists $x\in \limfunc{sol}(\mathrm{RP}_{\overline{x}^{\ast }})$
such that $J\in S_{f}\left( x\right) $ and $x_{j}^{\ast }\in \partial
f_{j}\left( x\right) $ for all $j\in J, $\newline
$(iv)$ There exists $x\in X$ such that $J\in S_{f}\left( x\right) $ and $%
x_{j}^{\ast }\in \partial f_{j}\left( x\right) $ for all $j\in J.$\newline
Moreover, for any $\left( J,(x_{j}^{\ast })_{j\in J}\right) \in \limfunc{sol}%
(\mathrm{RD}_{\overline{x}^{\ast }})$ one has
\begin{equation*}
\limfunc{sol}(\mathrm{RD}_{\overline{x}^{\ast }})=\left\{ \left(
J,(x_{j}^{\ast })_{j\in J}\right) \in \mathbb{F}\left( \overline{x}^{\ast
}\right) :J\in S_{f}\left( x\right) \text{ and }x_{j}^{\ast }\in \partial
f_{j}\left( x\right) ,\forall j\in J\right\}
\end{equation*}%
for either some (all) $x\in \limfunc{sol}(\mathrm{RP}_{\overline{x}^{\ast
}}) $ or for some $x\in X.$
\end{corollary}

\textit{Proof.} $\left[ (i)\Longrightarrow (ii)\right] $ It comes from the
statement $\left[ (i)\Longrightarrow (ii)\right] $ in Theorem \ref{Theor3.1}.

$\left[ (ii)\Longrightarrow (iii)\right] $ It is obvious as $\limfunc{sol}(%
\mathrm{RP}_{\overline{x}^{\ast }})\neq \emptyset .$

$\left[ (iii)\Longrightarrow (iv)\right] $ It is obvious.

$\left[ (iv)\Longrightarrow (i)\right] $ Since $\left( J,(x_{j}^{\ast
})_{j\in J}\right) \in \mathbb{F}\left( \overline{x}^{\ast }\right) ,$ $%
x_{j}^{\ast }\in \partial f_{j}\left( x\right) $ for all $j\in J,$ $%
\sum_{j\in J}x_{j}^{\ast }=\overline{x}^{\ast },$ and $J\in S_{f}\left(
x\right) ,$ one has%
\begin{equation*}
\begin{array}{ll}
\sup (\mathrm{RD}_{\overline{x}^{\ast }}) & \geq -\sum\limits_{j\in
J}f_{j}^{\ast }(x_{j}^{\ast }) \\
& =\sum\limits_{j\in J}\left( {{f_{j}}}\left( x\right) -\left\langle
x_{j}^{\ast },x\right\rangle \right) \\
& =\sum\limits_{j\in J}{{f_{j}}}\left( x\right) -\left\langle \overline{x}%
^{\ast },x\right\rangle \\
& =f\left( x\right) -\left\langle \overline{x}^{\ast },x\right\rangle \\
& \geq \inf (\mathrm{RP}_{\overline{x}^{\ast }}) \\
& \geq \sup (\mathrm{RD}_{\overline{x}^{\ast }}).%
\end{array}%
\end{equation*}%
Consequently, $\sup (\mathrm{RD}_{\overline{x}^{\ast }})=-\sum_{j\in
J}f_{j}^{\ast }(x_{j}^{\ast })$ and $(i)$ holds.

The last assertion of Corollary \ref{Corol3.2} comes directly from the
equivalences \noindent $(i)\Leftrightarrow (ii)\Leftrightarrow
(iii)\Leftrightarrow (iv).\hfill \square \medskip $

For the last result of this section we still assume $\left( f_{i}\right)
_{i\in I}\subset \left( \mathbb{R\cup }\left\{ +\infty \right\} \right) ^{X}$
is an infinite family of proper functions, but we do not consider a fixed
element $\overline{x}^{\ast }\in X^{\ast }.$ The equation (\ref{3.5}) is
called \textit{stable strong duality} in \cite{BGW09}.

\begin{corollary}
\label{Corol3.3}Assume that
\begin{equation}
\inf (\mathrm{RP}_{x^{\ast }})=\max (\mathrm{RD}_{x^{\ast }}),\forall
x^{\ast }\in \bigcup\limits_{x\in X}\partial f\left( x\right) .  \label{3.5}
\end{equation}%
Then one has%
\begin{equation}
\partial f\left( x\right) =\bigcup\limits_{J\in S_{f}\left( x\right)
}\sum\limits_{j\in J}\partial {{f_{j}}}\left( x\right) ,\forall x\in X.
\label{3.6}
\end{equation}
\end{corollary}

\textit{Proof.} Let us show that the inclusion $\supset $ always holds in (%
\ref{3.6}).

Let $x^{\ast }:=\sum_{j\in J}x_{j}^{\ast }$ with $J\in S_{f}\left( x\right) $
and $x_{j}^{\ast }\in \partial f_{j}\left( x\right) $ for all $j\in J.$ We
thus have,
\begin{equation*}
\begin{array}{ll}
f\left( x\right) -\left\langle x^{\ast },x\right\rangle & =\sum\limits_{j\in
J}\left( {{f_{j}}}\left( x\right) -\left\langle x_{j}^{\ast },x\right\rangle
\right) \\
& =-\sum\limits_{j\in J}f_{j}^{\ast }(x_{j}^{\ast }) \\
& \leq \sup (\mathrm{RD}_{x^{\ast }}) \\
& \leq \inf (\mathrm{RP}_{x^{\ast }}) \\
& =-f^{\ast }\left( x^{\ast }\right) \\
& \leq f\left( x\right) -\left\langle x^{\ast },x\right\rangle .%
\end{array}%
\end{equation*}%
Finally, $f\left( x\right) -\left\langle x^{\ast },x\right\rangle =-f^{\ast
}\left( x^{\ast }\right) ,$ that means $x^{\ast }\in \partial f\left(
x\right) .$

We now prove the reverse inclusion $\subset $ in (\ref{3.6}).

Let $x^{\ast }\in \partial f\left( x\right) .$ Then $x\in \partial f^{\ast
}\left( x^{\ast }\right) $ and, by (\ref{solRPx*}), $x\in \limfunc{sol}(%
\mathrm{RP}_{\overline{x}^{\ast }}).$ By (\ref{3.5}) and Corollary \ref%
{Corol3.1}, there exists $\left( J,(x_{j}^{\ast })_{j\in J}\right) \in
\mathbb{F}\left( \overline{x}^{\ast }\right) $ such that $J\in S_{f}\left(
x\right) $ and $x_{j}^{\ast }\in \partial f_{j}\left( x\right) $ for all $%
j\in J.$ We thus have $x^{\ast }=\sum\limits_{j\in J}x_{j}^{\ast }\in
\sum\limits_{j\in J}\partial {{f_{j}}}\left( x\right) . \hfill \square
\medskip $

\section{Robust sum of subaffine functions}

Let $\left( A_{i}\right) _{i\in I}\ $be a family of nonempty, $w^{\ast }-$%
closed convex subsets of $X^{\ast },$ $t_{i}\in \mathbb{R}$ for all $i\in I$
and the subaffine functions $f_{i}:=\sigma _{A_{i}}-t_{i},$ $i\in I.$ Then $%
(f_i)_{i \in I} \subset \Gamma (X) $ and we have $f_{i}^{\ast }:=\delta
_{A_{i}}+t_{i}$ and $\limfunc{epi}f_{i}^{\ast }=A_{i}\times \left[
t_{i},+\infty \right[ =A_{i}\times \left\{ t_{i}\right\} +\left\{ 0_{X^{\ast
}}\right\} \times \mathbb{R}_{+}$ for each $i \in I$. The robust sum $f$ of
this family is
\begin{equation*}
f\left( x\right) =\sum\nolimits_{i\in I}^{R}{{f_{i}}\left( x\right) =}%
\sup\limits_{J\in \mathcal{F}\left( I\right) }\sum\limits_{j\in J}\Big[%
\sigma _{Aj}\left( x\right) -t_{j}\Big],\forall x\in X
\end{equation*}
and the set $\mathcal{A}$ defined by \eqref{eqA} now becomes
\begin{equation}  \label{1.2}
\mathcal{A}:=\left( \dbigcup\limits_{\QATOP{J\in \mathcal{F}\left( I\right)
}{\hfill }}\sum\nolimits_{j\in J}\Big[A_{j}\times \left\{ t_{j}\right\} \Big]%
\right) +\left\{ 0_{X^{\ast }} \right\} \times \mathbb{R}_{+} .
\end{equation}

Let us introduce the set-valued mapping
\begin{equation*}
\mathbb{A}:\mathcal{F}\left( I\right) \rightrightarrows X^{\ast }\ \text{
such\ that \ } \mathbb{A}\left( J\right) =\sum\limits_{j\in J}A_{j}.
\end{equation*}
Then the problem $(\mathrm{RP}_{\overline{x}^*})$ and its dual $(\mathrm{RD}%
_{\overline{x}^*}) $ write as
\begin{equation*}
\inf (\mathrm{RP}_{\overline{x}^*}) = \inf \{ f(x) - \langle \overline{x}^*,
x \rangle : x \in X \} = - f^\ast (\overline{x}^*)
\end{equation*}
and
\begin{equation*}
\sup (\mathrm{RD}_{\overline{x}^*}) = \sup \left \{ - \sum_{j \in J}
f_i^\ast (x_j^*) : J \in \mathbb{A}^{-1} (\overline{x}^*) \right\} = - \inf
\left\{ \sum\limits_{j\in J}t_{j}: J\in \mathbb{A}^{-1} (\overline{x}^*)
\right\},
\end{equation*}
and hence, the zero duality gap relation amounts to%
\begin{equation*}
f^{\ast }\left( \overline{x}^{\ast }\right) =\inf \left\{ \sum_{j\in
J}t_{j}:J\in \mathbb{A}^{-1}\left( \overline{x}^{\ast }\right) \right\} .
\end{equation*}

We now briefly quote some remarkable properties on the duality and the
convexity and closedness of the qualifying set $\mathcal{A}$:

$\bullet$ It is worth observing firstly that if $\mathbb{A}^{-1}\left(
\overline{x}^{\ast }\right) =\emptyset$ (i.e., $\overline{x}^{\ast }\notin
\dbigcup_{ J\in \mathcal{F}\left( I\right)}\sum\limits_{j\in J}A_{j}$), one
has $\overline{x}^{\ast }\notin \limfunc{dom} f^\ast$ and $\sup (\mathrm{RD}%
_{\overline{x}^{\ast }})=-\infty .$

$\bullet $ In the case when $\func{dom}f\neq \emptyset $ (for instance, if $%
\sum\nolimits_{i\in I}^{R}t_{i}\in \mathbb{R}$), Theorem \ref{Theor1.1} says
that the \textit{stable strong duality of the pair $(\mathrm{RP}_{{x}^{\ast
}})$-$(\mathrm{RD}_{{x}^{\ast }})$ holds}, i.e.,
\begin{equation}
f^{\ast }\left( x^{\ast }\right) =\min \left\{ \sum\nolimits_{j\in
J}t_{j}:J\in \mathbb{A}^{-1}(x^{\ast })\right\} ,\forall x^{\ast }\in \func{%
dom}f^{\ast }  \label{1.3}
\end{equation}%
if and only if the set
\begin{equation}
\mathcal{A}=\left( \dbigcup\limits_{\QATOP{J\in \mathcal{F}\left( I\right) }{%
\hfill }}\sum\nolimits_{j\in J}\Big[A_{j}\times \left\{ t_{j}\right\} \Big]%
\right) +\left\{ 0_{X^{\ast }}\right\} \times \mathbb{R}_{+}\text{ is }\
w^{\ast }-\text{closed and convex}.  \label{closedness}
\end{equation}

$\bullet $ According to Lemma \ref{Lemma1.1} and Example \ref{Ex1.1}, we
know that the set $\mathcal{A}$ in (\ref{1.2}) is convex if $0_{X^{\ast
}}\in $ $\bigcap\limits_{i\in I}A_{i}\neq \emptyset $ and $%
\sup\nolimits_{i\in I}t_{i}\leq 0.$ Moreover, the set $\mathcal{A}$ is $%
w^{\ast }-$closed if $\dbigcup\limits_{\QATOP{J\in \mathcal{F}\left(
I\right) }{\hfill }}\sum\nolimits_{j\in J}\left( A_{j}\times \left\{
t_{j}\right\} \right) $ is $w^{\ast }-$compact.

On the primal attainment and the strong duality of the robust sum for
subaffine functions $(\mathrm{RP}_{x^{\ast }})$, one has the following
consequence of Theorem \ref{Theor2.1} and Lemma \ref{Lemma1.1}.

\begin{proposition}
\label{prop51nww} Assume that $0_{X^{\ast }}\in $ $\bigcap\limits_{i\in
I}A_{i}$ and the robust sum $\sum_{i\in I}^{R}\left( \sigma
_{A_{i}}-t_{i}\right) $ is proper and weakly inf-locally compact. Let $%
\overline{x}^{\ast }\in X^{\ast }$ be such that
\begin{equation}
\overline{\mathrm{cone}}\left( \bigcup\limits_{J\in \mathcal{F}%
(I)}\sum\limits_{j\in J}A_{j}-\overline{x}^{\ast }\right) \ \ \mathrm{is\ a\
linear\ subspace\ of}\ X^{\ast }.  \label{alpha}
\end{equation}%
Then the optimal solution set of the problem
\begin{equation*}
(\mathrm{RP}_{\overline{x}^{\ast }})\ \ \ \ \ \ \ \inf\limits_{x\in X}\left(
\sum^{R}\nolimits_{i\in I}(\sigma _{A_{i}}(x)-t_{i})-\left\langle \overline{x%
}^{\ast },x\right\rangle \right)
\end{equation*}%
is the sum of a nonempty weakly compact set and a finitely dimensional
linear subspace of $X$.
\end{proposition}

Applying Theorem \ref{Theor2.2} we get

\begin{proposition}
\label{prop52nww} Assume that $0_{X^{\ast }}\in $ $\bigcap\limits_{i\in
I}A_{i}\neq \emptyset $ and $\sup\nolimits_{i\in I}t_{i}\leq 0,$ and there
exists $i_{0}\in I$ such that $\delta _{A_{i_{0}}}$ is Mackey
quasicontinuous. Then for each $\overline{x}^{\ast }\in X^{\ast }$
satisfying \eqref{alpha} we have
\begin{equation*}
\min\limits_{x\in X}\left( \sum^{R}\nolimits_{i\in I}(\sigma
_{A_{i}}(x)-t_{i})-\left\langle \overline{x}^{\ast },x\right\rangle \right)
=\sup \big\{-\sum\limits_{j\in J}t_{j}:J\in \mathbb{A}^{-1}(\overline{x}%
^{\ast })\big\}.
\end{equation*}
\end{proposition}

\textit{Proof}. By Lemma \ref{Lemma1.1} (Example \ref{Ex1.1}) the set $%
\mathcal{A}$ is convex and the function $\varphi $ is convex, too (Remark %
\ref{rem32a}). On the other hand, by Remark \ref{rem32}, the function $%
\varphi $ is Mackey quasicontinuous. The conclusion follows from Theorem \ref%
{Theor2.2}. $\hfill \square \medskip $

In finite dimension we have (as an immediate consequence of Proposition \ref%
{prop52nww}):

\begin{proposition}
\label{prop53nww} Let $(A_{i})_{i\in I}$ be a family of closed convex
subsets of $\mathbb{R}^{n}$ such that $0_{X^{\ast }}\in $ $%
\bigcap\limits_{i\in I}A_{i}$. Assume that $\sup\nolimits_{i\in I}t_{i}\leq 0
$. Then for any $\overline{x}^{\ast }\in \limfunc{ri}\Big(\dbigcup\limits_{%
\QATOP{J\in \mathcal{F}\left( I\right) }{\hfill }}\sum\nolimits_{j\in J}A_{j}%
\Big)$ one has
\begin{equation*}
\min\limits_{x\in X}\left( \sum^{R}\nolimits_{i\in I}(\sigma
_{A_{i}}(x)-t_{i})-\left\langle \overline{x}^{\ast },x\right\rangle \right)
=\sup \big\{-\sum\limits_{j\in J}t_{j}:J\in \mathbb{A}^{-1}(\overline{x}%
^{\ast })\big\}.
\end{equation*}
\end{proposition}

We end this section with a formula on the subdifferential of the robust sum $%
f = \sum_{i\in I}^{R} (\sigma_{A_i} - t_i) $. Let us recall that for each $x
\in X$ one has, by definition,
\begin{equation*}
S_f (x) = \{ J \in \mathcal{F}(I) : \sum\limits_{j \in J} (\sigma_{A_j} (x)
- t_j) = f(x) \}.
\end{equation*}
We observe also that
\begin{equation*}
\partial \sigma_{A_i} (x) = \{ x^* \in A_i: \langle x^*, x \rangle =
\sigma_{A_i} (x) \}
\end{equation*}
or, in other words,
\begin{equation}  \label{eqnew}
\partial \sigma_{A_i} (x) = \mathrm{argmax}_{A_i} \langle \cdot , x \rangle.
\end{equation}
We then have:

\begin{proposition}
\label{prop54nww} Assume that $0_{X^{\ast }}\in $ $\bigcap\limits_{i\in
I}A_{i}$, $\sup\nolimits_{i\in I}t_{i}\leq 0$, $f$ is proper, and the set
\begin{equation}
\left( \dbigcup\limits_{\QATOP{J\in \mathcal{F}\left( I\right) }{\hfill }%
}\sum\nolimits_{j\in J}\Big[A_{j}\times \left\{ t_{j}\right\} \Big]\right)
+\left\{ 0_{X^{\ast }}\right\} \times \mathbb{R}_{+},  \label{eqneww}
\end{equation}%
is $w^{\ast }$-closed regarding the set $\bigcup\limits_{u\in X}\partial f(u)
$. Then one has
\begin{equation*}
\partial f(x)=\bigcup\limits_{J\in S_{f}(x)}\sum\limits_{j\in J}\mathrm{%
argmax}_{A_{j}}\langle \cdot ,x\rangle ,\forall x\in X.
\end{equation*}
\end{proposition}

\textit{Proof}. Noting that the set in \eqref{eqneww} is nothing but
\begin{equation*}
\mathcal{A} = \bigcup\limits_{J \in \mathcal{F}(I)} \sum\limits_{j \in J}
\limfunc{epi} (\sigma_{A_j} - t_j)^\ast ,
\end{equation*}
which is convex. The conclusion now follows from Theorem \ref{Theor1.1},
Corollary \ref{Corol3.3}, and \eqref{eqnew}. $\hfill \square \medskip $

\section{Approximate solutions to inconsistent convex inequality systems}

In this section $(f_{i})_{i\in I}\subset \Gamma (X)$. We consider the system
\begin{equation*}
\mathrm{\left( \mathrm{S}\right) }\hskip1cm\{f_{i}(x)\leq 0,\ i\in I\},
\end{equation*}%
that we assume to be inconsistent. Defining
\begin{equation}
f_{0}(x):=\sup\limits_{i\in I}f_{i}(x),  \label{f0}
\end{equation}%
we have $f_{0}(x)>0$ for all $x\in X$.

The $i-$th residual of $x$ is given by $f_{i}^{+}(x)$ and, in some sense,
the infeasibility of $x$ is measured by $\sup\limits_{i\in I}f_{i}^{+}(x),$
that is, $f_{0}(x)$ too. We may also consider the cumulative infeasibility
of $x$, namely the infinite sum $\sum\limits_{i\in I}f_{i}^{+}(x)$ (see \cite%
{GHL18}). Since $f_{0}(x)>0$ we know that $\sum\nolimits_{i\in I}f_{i}^{+}$
coincides with the robust sum $\sum_{i\in I}^{R}f_{i}$ of the family $%
(f_{i})_{i\in I}$ (see \cite[Lemma 2.5]{DGV19}).

In formal terms, let us define a best $\ell _{\infty }$-approximate solution
of the inconsistent system $\left( \mathrm{S}\right) $ as an optimal
solution to the problem
\begin{equation*}
\inf\nolimits_{x\in X}f_{0}(x)=\sup\limits_{i\in
I}f_{i}(x)=\sup\limits_{i\in I}f_{i}^{+}(x)
\end{equation*}%
and, similarly, a best $\ell _{1}$-approximate solution of $\left( \mathrm{S}%
\right) $ as an optimal solution to the problem
\begin{equation*}
\inf\nolimits_{x\in X}\sum\limits_{i\in I}f_{i}^{+}(x)=\sum\limits_{i\in
I}\nolimits^{R}f_{i}(x).
\end{equation*}%
We denote by $\ell _{\infty }$-$\limfunc{sol}\left( \mathrm{S}\right) $
(resp., $\ell _{1}$-$\limfunc{sol}\left( \mathrm{S}\right) $) the set of
best $\ell _{\infty }$ (resp., $\ell _{1}$) approximate solutions of the
inconsistent system $\left( \mathrm{S}\right) $.

In order to associate a suitable dual problem with $\inf\nolimits_{x\in
X}f_{0}(x)$ we define, as in \cite{GLV17}, the \textit{unit simplex} in the
linear space $\mathbb{R}^{\left( I\right) }$ of real-valued functions $%
\lambda \in \mathbb{R}^{I}$ with finite \textit{support set} $\limfunc{supp}%
\lambda :=\left\{ i\in I:\lambda _{i}\neq 0\right\} $ as
\begin{equation*}
S_{I}:=\left\{ \lambda \in \mathbb{R}^{\left( I\right) }:\sum\limits_{i\in
I}\lambda _{i}=1,\lambda _{i}\geq 0,\forall i\in I\right\}
\end{equation*}
and the \textit{modified Lagrangian function} as $L:X\times S_{I}$ such that
\begin{equation*}
L\left( x,\lambda \right) :=\sum\limits_{i\in \limfunc{supp}\lambda }\lambda
_{i}f_{i}(x),\forall \left( x,\lambda \right) \in X\times S_{I}.
\end{equation*}

\begin{proposition}[Structure of $\ell _{\infty }$-$\limfunc{sol}\left(
\mathrm{S}\right) $ and strong duality]
\label{Prop 6.1}Assume that $f_{0}$ is proper and weakly inf-locally
compact, and that $\overline{\limfunc{cone}}\limfunc{co}\bigcup\limits_{i\in
I}\limfunc{dom}f_{i}^{\ast }$ is a linear subspace of $X^{\ast }$. Then $%
\ell _{\infty }$ $-\limfunc{sol}\left( \mathrm{S}\right) $ is the sum of a
nonempty convex weakly compact set and a finitely dimensional linear
subspace of $X$. Moreover, one has%
\begin{equation*}
\inf\nolimits_{x\in X}\sup\limits_{i\in I}f_{i}(x)=\max \left\{ \inf_{x\in
X}\sum\limits_{i\in I}\lambda _{i}f_{i}(x):\lambda \in S_{I}\right\}
\end{equation*}%
if and only if
\begin{equation*}
\bigcup\limits_{\lambda \in S_{I}}\limfunc{epi}\left( \sum\limits_{i\in
I}\lambda _{i}f_{i}\right) ^{\ast }\text{ is }w^{\ast }-\text{closed
regarding }\{0_{X^{\ast }}\}\times \mathbb{R}.
\end{equation*}
\end{proposition}

\textit{Proof}. \ Since $f_{0}\in \Gamma (X)$ one has $\ell _{\infty }-%
\limfunc{sol}\left( \mathrm{S}\right) =\partial f_{0}^{\ast }(0_{X^{\ast }})$%
. We intend to apply Lemma \ref{Lemma2.1} for $g=f_{0}$ and $x^{\ast
}=0_{X^{\ast }}$. We have to make explicit the criterion \eqref{2.2} in
terms of the conjugate of the data functions $f_{i}$. To this end consider
the function $\Psi :=\inf\limits_{i\in I}f_{i}^{\ast }$. One has $\func{dom}%
\Psi =\cup _{i\in I}\func{dom}f_{i}^{\ast }$, $\Psi ^{\ast }=f_{0}$ and,
since $\func{dom}f_{0}\not=\emptyset $, $f_{0}^{\ast }=\overline{\limfunc{co}%
}\,\Psi $. Now, as in \eqref{2.5}, we have $\overline{\limfunc{co}}\func{dom}%
\Psi =\overline{\func{dom}\overline{\limfunc{co}}\Psi }$ and, consequently,
\begin{equation*}
\overline{\limfunc{cone}}\func{dom}f_{0}^{\ast }=\overline{\limfunc{cone}}({%
\func{dom}(\overline{\limfunc{co}}\Psi }))=\overline{\limfunc{cone}}\limfunc{%
co}(\func{dom}\Psi )=\overline{\limfunc{cone}}\Big(\limfunc{co}%
\bigcup\limits_{i\in I}\limfunc{dom}f_{i}^{\ast }\big).
\end{equation*}%
The strong duality theorem is consequence of \cite[Corollary 3.4]{GLV17}%
.\noindent $\hfill \square \medskip $

Observe that, if at least one of the\ functions $f_{i}$ is weakly
inf-locally compact, then $f_{0}$ is weakly inf-locally compact, too. The
next corollary is an immediate consequence of Proposition \ref{Prop 6.1}.

\begin{corollary}
\label{Corol6.1}Assume that $(f_{i})_{i\in I}\subset \Gamma (\mathbb{R}%
^{n}), $ $\limfunc{dom}f_{0}\neq \emptyset $, and $0_{\mathbb{R}^{n}}\in
\limfunc{ri}\limfunc{co}\Big(\bigcup\limits_{i\in I}\limfunc{dom}f_{i}^{\ast
}\Big).$ Then $\ell _{\infty }$-$\limfunc{sol}\left( \mathrm{S}\right) $ is
the sum of a nonempty convex compact set and a linear subspace of $\mathbb{R}%
^{n}$.
\end{corollary}

\begin{example}
\label{Ex6.1}Let $\left\{ \left\langle a_{i},x\right\rangle \leq b_{i},i\in
I\right\} $ be an inconsistent linear system posed in $\mathbb{R}^{n}.$ This
is a particular case of system $\left( \mathrm{S}\right) $\ above, with $%
f_{i}=\left\langle a_{i},\cdot \right\rangle -b_{i},$ $a_{i}\in \mathbb{R}%
^{n}$ and $b_{i}\in \mathbb{R}$ for all $i\in I.$ Denoting by $0_{n}$\ the
null vector in $\mathbb{R}^{n},$\ by Corollary \ref{Corol6.1}, if $\limfunc{%
dom}f_{0}\neq \emptyset $ and $0_{n}\in \limfunc{ri}\limfunc{co}\left\{
a_{i},i\in I\right\} ,$ then $\ell _{\infty }$-$\limfunc{sol}\left( \mathrm{S%
}\right) $ is the sum of a nonempty convex compact set and a linear subspace
of $\mathbb{R}^{n}$ (\cite[Proposition 1(S)]{GHL18} only asserts that, under
these assumptions, $\ell _{\infty }$-$\limfunc{sol}\left( \mathrm{S}\right)
\neq \emptyset $). Moreover, since
\begin{equation*}
\bigcup\limits_{\lambda \in S_{I}}\limfunc{epi}\left( \sum\limits_{i\in
I}\lambda _{i}f_{i}\right) ^{\ast }=\left\{ \sum\limits_{i\in I}\lambda
_{i}\left( a_{i},b_{i}\right) :\lambda \in S_{I}\right\} +\{0_{n}\}\times
\mathbb{R}_{+},
\end{equation*}%
the strong duality theorem becomes here
\begin{equation*}
\inf\nolimits_{x\in \mathbb{R}^{n}}\sup\limits_{i\in I}\left( \left\langle
a_{i},x\right\rangle -b_{i}\right) =\max \left\{ \inf_{x\in \mathbb{R}%
^{n}}\sum\limits_{i\in I}\lambda _{i}\left( \left\langle
a_{i},x\right\rangle -b_{i}\right) :\lambda \in S_{I}\right\} ,
\end{equation*}%
if and only if
\begin{equation*}
\left\{ \sum\limits_{i\in I}\lambda _{i}\left( a_{i},b_{i}\right) :\lambda
\in S_{I}\right\} +\{0_{n}\}\times \mathbb{R}_{+}\text{ is closed regarding }%
\{0_{n}\}\times \mathbb{R}_{+}.
\end{equation*}
\end{example}

\begin{proposition}[Structure of $\ell _{1}$-$\limfunc{sol}\left( \mathrm{S}%
\right) $ and strong duality]
\label{Prop 6.2}Assume that the robust sum $\sum\nolimits_{i\in I}^{R}f_{i}$
is proper, weakly inf-locally compact, and $\overline{\limfunc{cone}}%
\limfunc{co}\Big(\bigcup\limits_{J\in \mathcal{F}(I)}\sum\limits_{j\in J}%
\limfunc{dom}f_{j}^{\ast }\big)\ ${is\ a\ linear\ subspace\ of\ }$X^{\ast }.$
Then, $\ell _{1}$-$\limfunc{sol}\left( \mathrm{S}\right) $ is the sum of a
nonempty convex weakly compact set and a finitely dimensional linear
subspace of $X$. Moreover, one has
\begin{equation*}
\inf\nolimits_{x\in X}\sum_{i\in I}f_{i}^{+}(x)=\max \Big\{-\sum_{j\in
J}f_{j}^{\ast }(x_{j}^{\ast }):J\in \mathcal{F}(I),(x_{j}^{\ast })_{j\in
J}\in (X^{\ast })^{J},\sum_{j\in J}x_{j}^{\ast }=0_{X^{\ast }}\Big\}
\end{equation*}%
if and only if%
\begin{equation*}
\bigcup_{J\in \mathcal{F}(I)}\sum_{j\in J}\limfunc{epi}f_{j}^{\ast }\text{
is }w^{\ast }-\text{closed convex regarding }\{0_{X^{\ast }}\}\times \mathbb{%
R}.
\end{equation*}
\end{proposition}

\textit{Proof}. \ It is direct consequence of Theorem \ref{Theor2.1} and
Theorem \ref{Theor1.1} for $x^{\ast }=0_{X^{\ast }}$, due to the relation $%
\sum\nolimits_{i\in I}^{R}f_{i}=\sum\nolimits_{i\in I}f_{i}^{+}$. \noindent $%
\hfill \square \medskip $

\begin{example}
Consider again the linear system $\left( \mathrm{S}\right) $ in Example \ref%
{Ex6.1}. By Proposition \ref{Prop 6.2}, if $\sum\nolimits_{i\in I}^{R}\left(
\left\langle a_{i},\cdot \right\rangle -b_{i}\right) $ is proper and $%
0_{n}\in \limfunc{ri}\Big(\bigcup\limits_{J\in \mathcal{F}%
(I)}\sum\limits_{j\in J}a_{j}\big),\ $then $\ell _{1}$-$\limfunc{sol}\left(
\mathrm{S}\right) $ is the sum of a nonempty convex compact set and a
finitely dimensional linear subspace of $\mathbb{R}^{n}$. Observe that, for
each $(x_{j})_{j\in J}\in \left( \mathbb{R}^{n}\right) ^{J},$ one has
\begin{equation*}
\sum_{j\in J}f_{j}^{\ast }(x_{j})=\sum_{j\in J}\left( \delta _{\left\{
a_{j}\right\} }^{\ast }(x_{j})+b_{j}\right) =\left\{
\begin{array}{ll}
\sum\limits_{j\in J}b_{j}, & \text{if }x_{j}=a_{j},\text{ }\forall j\in J,
\\
+\infty , & \text{else.}%
\end{array}%
\right.
\end{equation*}%
So, again by Proposition \ref{Prop 6.2},%
\begin{equation*}
\inf\nolimits_{x\in \mathbb{R}^{n}}\sum_{i\in I}\left( \left\langle
a_{i},x\right\rangle -b_{i}\right) ^{+}=\max \Big\{-\sum_{J\in \mathcal{F}%
(I)}b_{j}:J\in \mathcal{F}(I),\sum_{J\in \mathcal{F}(I)}a_{j}=0_{n}\Big\}
\end{equation*}%
if and only if%
\begin{equation*}
\bigcup_{J\in \mathcal{F}(I)}\sum_{j\in J}\left( \{a_{i}\}\times \left[
b_{i},+\infty \right[ \right) \text{ is closed convex regarding }%
\{0_{n}\}\times \mathbb{R}.
\end{equation*}
\end{example}

\textbf{Acknowledgements} The authors wish to thank two anonymous referees
for their valuable comments which helped to improve the manuscript. \newline
This research was supported by the National Foundation for Science \&
Technology Development (NAFOSTED), Vietnam, Project 101.01-2018.310 \textit{%
Some topics on systems with uncertainty and robust optimization}, and by the
Ministry of Science, Innovation and Universities of Spain and the European
Regional Development Fund (ERDF) of the European Commission, Project
PGC2018-097960-B-C22.

\end{document}